\input amstex

\documentstyle{amsppt}

\refstyle{A}

\nologo


\hoffset .35 true in
\voffset .2 true in

\hsize=5.9 true in
\vsize=8.3 true in

\topmatter

\title Integral Restrictions on the Monodromy \endtitle

\author David B. Massey \endauthor

\address{David B. Massey, Dept. of Mathematics, Northeastern University, Boston, MA, 02115, USA} \endaddress

\email{DMASSEY\@NEU.edu}\endemail


\abstract Given a complex analytic function with a one-dimensional critical locus at the origin, we examine the monodromy action on the integral cohomology of the Milnor fiber. We relate this monodromy to that of a generic hyperplane slice through the origin, and to that of a generic hyperplane slice near the origin. We thereby obtain number-theoretic restrictions on the monodromy and on the cohomology of the original Milnor fiber.   \endabstract

\endtopmatter

\document

\baselineskip= 14pt

\noindent\S1. {\bf Introduction}  

\vskip .1in

Let $\Cal U$ be an open neighborhood of the origin in $\Bbb C^{n+1}$, and let $z_0$ be a non-zero linear form on $\Bbb C^{n+1}$. Let $f:(\Cal U, \bold 0)\rightarrow(\Bbb C, 0)$ be a complex analytic function with a $1$-dimensional critical locus at the origin, i.e., ${\operatorname{dim}}_\bold 0\Sigma f = 1$. We assume that $\Cal U$ is chosen small enough so that $\Sigma f\subseteq V(f)$. 

Let $H_0:=\Cal U\cap V(z_0)$ and let $H_t:=\Cal U\cap V(z_0-t)$, where $t\neq 0$ is sufficiently small. Let $f_0:= f_{|_{H_0}}$ and $f_t:=f_{|_{H_t}}$. We assume that $z_0$ is generic enough so that $\Sigma(f_0) = \{\bold 0\}$. This implies that $\Sigma(f_t)$ consists of a finite number of points: the points of $H_t\cap\Sigma f$, and the points where $H_t$ intersects the relative polar curve $\Gamma^1_{f, z_0}$.

Let $F$ denote the Milnor fiber of $f$ at $\bold 0$, and let $m_{n-1}: \widetilde H^{n-1}(F; \ \Bbb Z)@>\cong>> \widetilde H^{n-1}(F; \ \Bbb Z)$ and $m_{n}: \widetilde H^{n}(F; \ \Bbb Z)@>\cong>> \widetilde H^{n}(F; \ \Bbb Z)$ denote the corresponding $f$-monodromy actions on reduced integral cohomology. Note that $\widetilde H^{n-1}(F; \ \Bbb Z)$ is free-Abelian. Let $F_0$ denote the Milnor fiber of $f_0$ at $\bold 0$, and let $$h_0: \Bbb Z^{\mu_0}\cong \widetilde H^{n-1}(F_0; \ \Bbb Z)@>\cong>> \widetilde H^{n-1}(F_0; \ \Bbb Z)\cong\Bbb Z^{\mu_0}$$  denote the corresponding $f_0$-monodromy action, where $\mu_0$ denotes the Milnor number of the isolated critical point. For each point $\bold p_i\in H_t\cap\Sigma f$, there is an associated Milnor fiber $F_i$ of $f_t$ at $\bold p_i$, together with the associated $f_t$-monodromy action on $$h_i: \Bbb Z^{\mu_i}\cong \widetilde H^{n-1}(F_i;\ \Bbb Z)@>\cong>> \widetilde H^{n-1}(F_i;\ \Bbb Z)\cong \Bbb Z^{\mu_i}.$$

The monodromies of $f$ restricted to $H_0$ and $H_t$ are compatible with the monodromy of $f$ itself; they are all determined by letting the value of $f$ move in a small circle around the origin in $\Bbb C$. In categorical terms, this compatibility is a result of the fact that the monodromy is a natural automorphism of the vanishing cycle functor $\phi_f$.

This compatibility enables us to produce a diagram on which the monodromy acts
$$
\CD
{}@. \bigoplus_i\Bbb Z^{\mu_i} @.{} \\
@. @VV\alpha V @. \\
\Bbb Z^{\mu_0}@>\gamma>>N@>\delta>>\Bbb Z^{\lambda^0_{f}}\\
@. @VV\beta V @. \\
{}@. \Bbb Z^{\omega}, @.{}
\endCD
$$
where both the row and the column are short exact (we have omitted the zeroes on each end), $\dsize\lambda^0_f= \left(\Gamma^1_{f, z_0}\cdot V\left(\frac{\partial f}{\partial z_0}\right)\right)_\bold 0$ is the $0$-th L\^e number of $f$ (with respect to $z_0$ at the origin), and $\omega:= \big(\Gamma^1_{f, z_0}\cdot V(f)\big)_\bold 0$. It is important, and well-known, that $\omega\geqslant\lambda^0_f$, with equality if and only if $\omega =\lambda^0_f = 0$.

It is also important to note that, by the result of A'Campo [{\bf A}], the trace of of the monodromy maps $h_0$ and $h_i$ are all equal to $(-1)^n$.

\vskip .1in

We refer to the $N$ in the middle of the diagram as the {\it nexus}, and refer to the whole diagram as the {\it nexus diagram}. In fact, it is fairly unimportant what the nexus of the diagram is; what is important is that there exists such a diagram on which the $f$-monodromy acts. 

Let us denote the pull-back via $\alpha$ and $\gamma$ by $P$, and the push-forward via $\beta$ and $\delta$ by $Q$. Then, it is trivial to show that 
$$
P\ \cong \  \operatorname{ker}(\beta\circ \gamma) \ \cong  \ \operatorname{ker}(\delta\circ \alpha)
$$
and
$$
Q\ \cong \  \operatorname{coker}(\beta\circ \gamma) \ \cong  \ \operatorname{coker}(\delta\circ \alpha)
$$
This is important because $\beta\circ\gamma$ is the map induced on cohomology by the Morse-theoretic attaching map given by L\^e in [{\bf L1}]. This implies that $\operatorname{ker}(\beta\circ \gamma)\cong \widetilde H^{n-1}(F;\  \Bbb Z)$ and $\operatorname{coker}(\beta\circ \gamma)\cong \widetilde H^{n}(F; \ \Bbb Z)$. An alternate way of seeing these isomorphisms is in terms of L\^e numbers (see [{\bf Ma2}]); the $\Bbb Z$-module $\bigoplus_i\Bbb Z^{\mu_i} $ is precisely $\Bbb Z^{\lambda^1_f}$, where $\lambda^1_f$ is the $1$-dimensional L\^e number. The map $\delta\circ\alpha:\Bbb Z^{\lambda^1_f}\rightarrow \Bbb Z^{\lambda^0_f}$ is precisely the integral version of the one non-trivial map appearing in the chain complex of Corollary 10.10 of [{\bf Ma2}] (though there is a typographically error in this corollary -- the arrows are reversed). As the cohomology of this complex yields the reduced cohomology of $F$, we again conclude that $P\cong \widetilde H^{n-1}(F;\  \Bbb Z)$ and $Q\cong \widetilde H^{n}(F;\  \Bbb Z)$.

\vskip .1in

However, it is not only the $f$-monodromy that acts on the nexus diagram. We shall see that the $z_0$-monodromy also acts commutatively on the diagram, acting as the identity on the $\Bbb Z^{\mu_0}$ and $\Bbb Z^\omega$ nodes. Moreover, since the $z_0$ monodromy is also a natural isomorphism, it follows that the actions of the $f$-monodromy and $z_0$-monodromy on the nexus diagram commute with each other.

To continue our discussion, we must adopt some notation for the $z_0$-monodromy action on $ \bigoplus_i\Bbb Z^{\mu_i}$. Let $\nu$ be a (reduced) component of $\Sigma f$. Assume that $\Cal U$ is small enough so that $\nu$ is homeomorphic to a disk, and small enough so that the vanishing cycles along $f$, restricted to the punctured disk $\nu^*:=\nu-\{\bold 0\}$, form a local system. Then, there is an {\it internal} (or {\it vertical}) monodromy action $\iota_\nu$ induced on a given stalk of the vanishing cycle local system; one moves once around the ``hole'' in the punctured disk $\nu^*$.

Now, for each component $\nu$ of $\Sigma f$, there are $k_\nu:=(\nu\cdot V(z_0))_\bold 0$ points $\bold p_i$ which occur in $\nu\cap H_t$. At each of these $\bold p_i$, $f_t$ has the same Milnor number; let us denote this common value by $\mu_\nu$. There is a ``fractional monodromy'' action, $\tau_\nu:\Bbb Z^{\mu_\nu}@>\cong>> \Bbb Z^{\mu_\nu}$ given by moving cyclicly from one $\bold p_i$  in $\nu\cap H_t$ to the next. It follows that $\tau_{\nu}^{k_\nu}=\iota_v$. Moreover, if we let $\lambda_\nu$ denote the $z_0$-monodromy action on $\dsize \bigoplus_{\bold p_i\in \nu\cap H_t}\Bbb Z^{\mu_i}$, then it follows that $\lambda_\nu:\left(\Bbb Z^{\mu_\nu}\right)^{k_\nu}\rightarrow \left(\Bbb Z^{\mu_\nu}\right)^{k_\nu}$ is given by
$$
\lambda_\nu(\bold v_1, \bold v_2, \dots, \bold v_{k_\nu}) = (\tau_{\nu}(\bold v_{k_\nu}), \tau_{\nu}(\bold v_1), \dots, \tau_{\nu}(\bold v_{k_\nu-1})),
$$
from which one concludes immediately that $\operatorname{ker}(\operatorname{id}-\iota_\nu)\cong \operatorname{ker}(\operatorname{id}-\lambda_\nu)$.

\vskip .1in

Returning to the nexus diagram, since the $z_0$-monodromy acts as the identity on $\Bbb Z^{\mu_0}$, it follows that $\widetilde H^{n-1}(F; \ \Bbb Z)$ must be contained in $\bigoplus_\nu \operatorname{ker}(\operatorname{id}-\lambda_\nu)\cong \bigoplus_\nu \operatorname{ker}(\operatorname{id}-\iota_\nu)$. In particular, the rank of $\widetilde H^{n-1}(F; \ \Bbb Z)$ is at most $\sum_\nu \mu_\nu$.

Let us use ${\operatorname{char}}_{{}_-}(t)$ to denote characteristic polynomials. Then, using our discussion above and the proof in Section 2, what we show is:

\vskip .3in

\noindent{\bf Main Theorem}. {\it The nexus diagram exists, and the $f$-monodromy and $z_0$-monodromy act commutatively on it, and these actions commute with each other. 

The $z_0$-mondromy acts as the identity on $\Bbb Z^{\mu_0}$ and on $\Bbb Z^\omega$. The kernel of the identity minus the $z_0$-monodromy on $\bigoplus_i\Bbb Z^{\mu_i}$ is isomorphic to $\bigoplus_\nu \operatorname{ker}(\operatorname{id}-\iota_\nu)$.

Therefore, $\widetilde H^{n-1}(F; \ \Bbb Z)$ injects into $\bigoplus_\nu \operatorname{ker}(\operatorname{id}-\iota_\nu)$. Moreover, if $h_\nu$ denotes one of the $h_i$ for $\bold p_i\in\nu$, then ${\operatorname{char}}_{m_{n-1}}(t)$ divides ${\operatorname{char}}_{h_0}(t)$ and $\prod_\nu{\operatorname{char}}_{h_\nu}(t)$ in $\Bbb Z[t]$, i.e.,
$$
{\operatorname{char}}_{m_{n-1}}(t)\ \big|\ \operatorname{gcd}\Big({\operatorname{char}}_{h_0}(t), \ \prod_\nu{\operatorname{char}}_{h_\nu}(t)\Big).
$$
}

\vskip .2in

Aside from proving the Main Theorem, we also give a few applications of it. We show how the nexus diagram simplifies the proof of the non-splitting result of L\^e in [{\bf L2}]. We show how the Main Theorem allows us to generalize L\^e's non-splitting result and obtain the main result of [{\bf Ma1}].  Finally, we show how the Main Theorem generalizes the main result of [{\bf Ma3}] from the case of arrangements of planes in $\Bbb C^3$ to the case of arbitrary affine hypersurfaces with $1$-dimensional critical loci.

In the final section of the paper, we make some brief remarks which may lead to future applications.

\vskip .3in

\noindent\S2. {\bf Proof of the Main Theorem}.

\vskip .1in

We continue with the notation from the introduction.

\vskip .1in

Let $k:\Cal U\cap V(z_0)\hookrightarrow\Cal U$, $\hat k:V(f)\cap V(z_0)\hookrightarrow V(f)$, and $l:\Cal U-\Cal U\cap V(z_0)\hookrightarrow\Cal U$ denote the inclusions. Let . Let $\hat z_0 := {z_0}_{|_{V(f)}}$. 

Let $\bold P^\bullet$ denote the complex of sheaves of $\Bbb Z$-modules $\Bbb Z^\bullet_\Cal U[n+1]$; this sheaf is perverse. There is a fundamental distinguished triangle
$$
l_!l^!\bold P^\bullet\rightarrow\bold P^\bullet\rightarrow k_*k^*\bold P^\bullet@>[1]>>,
$$
which we can ``turn'' to yield
$$
k_*k^*\bold P^\bullet[-1]\rightarrow l_!l^!\bold P^\bullet\rightarrow\bold P^\bullet@>[1]>>.
$$
Applying the composed functor $\phi_{\hat z_0}[-1]\phi_f[-1]$ to the above triangle, we obtain the distinguished triangle
$$
\phi_{\hat z_0}[-1]\phi_f[-1]k_*k^*\bold P^\bullet[-1]\rightarrow \phi_{\hat z_0}[-1]\phi_f[-1]l_!l^!\bold P^\bullet\rightarrow\phi_{\hat z_0}[-1]\phi_f[-1]\bold P^\bullet@>[1]>>.\tag{$\dagger$}
$$

\vskip .1in

Now, 

\vskip .1in

\noindent$\bullet$\hskip .2in $k^*\bold P^\bullet[-1]\cong \Bbb Z^\bullet_{\Cal U\cap V(z_0)}[n]$, and so $\phi_f[-1]k_*k^*\bold P^\bullet[-1]\cong \hat k_*\phi_{f_0}[-1]\Bbb Z^\bullet_{\Cal U\cap V(z_0)}[n]$. As the support of this last complex of sheaves is contained in $V(z_0)$, if we apply $\psi_{\hat z_0}[-1]$, we get the zero complex. Hence, we find that the first complex in $(\dagger)$ 
$$
\phi_{\hat z_0}[-1]\phi_f[-1]k_*k^*\bold P^\bullet[-1]\cong \phi_{f_0}[-1]\Bbb Z^\bullet_{\Cal U\cap V(z_0)}[n].
$$
Note that this is a perverse sheaf with isolated support at the origin.

\vskip .1in

\noindent$\bullet$\hskip .2in As $\bold P^\bullet$ is perverse and $l$ is the inclusion of a hypersurface complement, $l_!l^!\bold P^\bullet$ is perverse. Thus, the second complex of $(\dagger)$, $\phi_{\hat z_0}[-1]\phi_f[-1]l_!l^!\bold P^\bullet$,  is perverse.

\vskip .1in

\noindent$\bullet$\hskip .2in As $\bold P^\bullet$ is perverse and the origin is an isolated point in the intersection of $V(z_0)$ and $\Sigma f$, $\phi_{\hat z_0}[-1]\phi_f[-1]\bold P^\bullet$ is perverse and has the origin as an isolated point in its support.

\vskip .2in

At isolated points in their supports, perverse sheaves have their stalk cohomology concentrated in degree zero. Thus, the long-exact sequence on the stalk cohomology at the origin obtained from $(\dagger)$ has at most one non-trivial piece -- namely, the short exact sequence
$$
0\rightarrow H^0\big(\phi_{f_0}[-1]\Bbb Z^\bullet_{\Cal U\cap V(z_0)}[n]\big)_\bold 0\rightarrow H^0\big(\phi_{\hat z_0}[-1]\phi_f[-1]l_!l^!\bold P^\bullet\big)_\bold 0\rightarrow H^0\big(\phi_{\hat z_0}[-1]\phi_f[-1]\bold P^\bullet\big)_\bold 0\rightarrow 0.\tag{$\ddagger$}
$$

Observe that $H^0\big(\phi_{f_0}[-1]\Bbb Z^\bullet_{\Cal U\cap V(z_0)}[n]\big)_\bold 0\cong\Bbb Z^{\mu_0}$ and that, by [{\bf Ma2}] or [{\bf Ma4}], 
$$
 H^0\big(\phi_{\hat z_0}[-1]\phi_f[-1]\Bbb Z^\bullet_\Cal U[n+1]\big)_\bold 0\cong \Bbb Z^{\lambda^0_f}.
$$
Therefore, we define the nexus $N$ to be $H^0\big(\phi_{\hat z_0}[-1]\phi_f[-1]l_!l^!\bold P^\bullet\big)_\bold 0$, and $(\ddagger)$ becomes the horizontal row of the nexus diagram.

\vskip .2in

To obtain the vertical row of the nexus diagram, we begin by considering another fundamental distinguished triangle
$$
\hat k^*[-1]\phi_f[-1]l_!l^!\bold P^\bullet\rightarrow\psi_{\hat z_0}[-1]\phi_f[-1]l_!l^!\bold P^\bullet\rightarrow\phi_{\hat z_0}[-1]\phi_f[-1]l_!l^!\bold P^\bullet@>[1]>>
$$
and turn this triangle to obtain
$$
\psi_{\hat z_0}[-1]\phi_f[-1]l_!l^!\bold P^\bullet\rightarrow\phi_{\hat z_0}[-1]\phi_f[-1]l_!l^!\bold P^\bullet\rightarrow \hat k^*\phi_f[-1]l_!l^!\bold P^\bullet@>[1]>>.
$$

Again, we look at the associated long exact sequence on the stalk cohomology at the origin.

\vskip .1in

\noindent$\bullet$\hskip .2in As $l_!l^!\bold P^\bullet$ is isomorphic to $\bold P^\bullet$ outside of $V(z_0)$, it follows that $$\psi_{\hat z_0}[-1]\phi_f[-1]l_!l^!\bold P^\bullet\ \cong\ \psi_{\hat z_0}[-1]\phi_f[-1]\bold P^\bullet.$$  Since $\bold P^\bullet$ is perverse and the origin is an isolated point in the intersection of $V(z_0)$ and $\Sigma f$, $\psi_{\hat z_0}[-1]\phi_f[-1]\bold P^\bullet$ is perverse and has the origin as an isolated point in its support. Thus,
$$
H^0\big(\psi_{\hat z_0}[-1]\phi_f[-1]l_!l^!\bold P^\bullet\big)_\bold 0\ \cong\ H^0\big(\psi_{\hat z_0}[-1]\phi_f[-1]\bold P^\bullet\big)_\bold 0\ \cong\ \bigoplus_i\Bbb Z^{\mu_i}.
$$

\vskip .1in

\noindent$\bullet$\hskip .2in Note that $H^*\big(\hat k^*\phi_f[-1]l_!l^!\bold P^\bullet)_0\cong H^*\big(\phi_f[-1]l_!l^!\bold P^\bullet)_0$ As the nexus is defined to be $N=H^0\big(\phi_{\hat z_0}[-1]\phi_f[-1]l_!l^!\bold P^\bullet\big)_\bold 0$, we would be finished if we could show that $H^*\big(\phi_f[-1]l_!l^!\bold P^\bullet)_0$ is zero outside of degree zero and that 
$$H^0\big(\phi_f[-1]l_!l^!\bold P^\bullet)_0\cong \Bbb Z^\omega.$$

Consider the distinguished triangle
$$
\big(l_!l^!\bold P^\bullet\big)_{|_{V(f)}}\rightarrow \psi_f[-1]l_!l^!\bold P^\bullet\rightarrow\phi_f[-1]l_!l^!\bold P^\bullet@>[1]>>.
$$
As $l_!$ is the extension by zero, we find that
$$H^0\big(\phi_f[-1]l_!l^!\bold P^\bullet)_0\cong H^0\big(\psi_f[-1]l_!l^!\bold P^\bullet)_0.$$

Now, consider the distinguished triangle
$$
\psi_f[-1]l_!l^!\bold P^\bullet\rightarrow\psi_f[-1]\bold P^\bullet@>\xi>> \psi_f[-1]k_*k^*\bold P^\bullet@>[1]>>.
$$
The stalk at $\bold 0$ of the map $\xi$ is the map induced by L\^e's attaching map in [{\bf L1}]; this is Morse-theoretic map involved in L\^e's description of how the Milnor fiber of $f$ is obtained from the Milnor fiber of $f_0$. Therefore, 
$H^*\big(\psi_f[-1]l_!l^!\bold P^\bullet)_0$ is isomorphic to the relative cohomology $H^*(F, F_0;\ \Bbb Z)$, and L\^e's result tells us that $H^i\big(\psi_f[-1]l_!l^!\bold P^\bullet)_0 =0$, unless $i=0$, and
$$H^0\big(\psi_f[-1]l_!l^!\bold P^\bullet)_0\cong \Bbb Z^\omega.$$
This is what we needed to prove.

\vskip .3in

\noindent\S3. {\bf Applications}.

\vskip .1in

\noindent{\bf Application 1}. We will use the Main Theorem to simplify the non-splitting result of L\^e which appears in [{\bf L2}]. The heart of the argument is same; it is the details and ``trick'' at the end of L\^e's proof that we can eliminate.

\vskip .1in

Suppose that $\mu_0 = \lambda^1_f$, i.e., that $\mu_0 = \sum_i\mu_i$. We wish to conclude that there is exactly one point in $H_t\cap \Sigma f$, i.e., that the intersection number of reduced varieties $\big(|\Sigma f|\cdot V(z_0)\big)_\bold 0$ equals one. This is equivalent to saying that $\Sigma f$ has a single smooth component which is transversely intersected by $V(z_0)$.

\vskip .1in

As $\mu_0 = \lambda^1_f$, from the nexus diagram we conclude that $\omega = \lambda^0_f$. However, as discussed in the Introduction, this implies that $\omega = \lambda^0_f = 0$. Thus, $\alpha$ and $\gamma$ are isomorphisms. As the $f$-monodromy maps act commutatively on the nexus diagram, the result of A'Campo implies that 
$$
(-1)^n = \big(|\Sigma f|\cdot V(z_0)\big)_\bold 0(-1)^n,
$$
and hence $\big(|\Sigma f|\cdot V(z_0)\big)_\bold 0=1$.

\vskip .1in

Note that the nexus diagram also implies that, in this case, $\widetilde H^{n}(F; \ \Bbb Z) = 0$ and $\widetilde H^{n-1}(F; \ \Bbb Z)\cong \widetilde H^{n-1}(F_0; \ \Bbb Z)$.

\vskip .3in

\noindent{\bf Application 2}. The fact that $\widetilde H^{n-1}(F; \ \Bbb Z)\cong\operatorname{ker}(\delta\circ \alpha)$ immediately implies that the rank of $\widetilde H^{n-1}(F; \ \Bbb Z)$ is at most $\lambda^1_f$. Assume that $\operatorname{rk}\widetilde H^{n-1}(F; \ \Bbb Z)= \lambda^1_f$; we wish to see what this implies.

\vskip .1in

 Recall from the Introduction that $\operatorname{rk}\widetilde H^{n-1}(F; \ \Bbb Z)\leqslant\sum_\nu \mu_\nu$. As $\lambda^1_f= \sum_\nu (\nu\cdot V(z_0))_\bold 0\,\mu_{\nu}$, we see that we must have that, for every component $\nu$ of $\Sigma f$, $(\nu\cdot V(z_0))_\bold 0 =1$, i.e., each $\nu$ is smooth at the origin and is transversely intersected by $V(z_0)$.

 Now, $\widetilde H^{n-1}(F; \ \Bbb Z)\cong \Bbb Z^{\lambda^1_f}$ injects into the copy of $\Bbb Z^{\lambda^1_f}$ in the nexus diagram. Hence, the characteristic polynomials of the $f$-monodromy on $\widetilde H^{n-1}(F; \ \Bbb Z)$ and on the copy of $\Bbb Z^{\lambda^1_f}$ in the nexus diagram must be equal. Using that $\widetilde H^{n-1}(F; \ \Bbb Z)$ injects into  $\widetilde H^{n-1}(F_0; \ \Bbb Z)$ and applying the result of A'Campo, we conclude that $\widetilde H^{n-1}(F_0; \ \Bbb Z)\cong \Bbb Z^{\mu_0}$ has an $f$-monodromy-invariant free submodule of rank $\lambda^1_f$ on which the $f$-monodromy acts with trace $\big(|\Sigma f|\cdot V(z_0)\big)_\bold 0(-1)^n$. Since the eigenvalues of the $f$-monodromy are all roots of unity, and the trace of $h_0$ is itself equal to $(-1)^n$, we find that
$$
\big(|\Sigma f|\cdot V(z_0)\big)_\bold 0(-1)^n +\big(\text{the sum of }\mu_0-\lambda^1_f\text{ roots of unity}\big) = (-1)^n,
$$
or 
$$
\Big(\big(|\Sigma f|\cdot V(z_0)\big)_\bold 0-1\Big)(-1)^n \ =\ \text{the sum of }\mu_0-\lambda^1_f\text{ roots of unity}.\tag{$\dagger$}
$$
Taking norms, we conclude that 
$$
\big(|\Sigma f|\cdot V(z_0)\big)_\bold 0-1\leqslant \mu_0-\lambda^1_f\tag{$\ddagger$}
$$
with equality implying that $h_0$ has $\mu_0-\lambda^1_f$ eigenvalues equal to $(-1)^{n+1}$. 

\vskip .2in

\noindent$\bullet$\hskip .2in The case in Application 1 was that of $\mu_0-\lambda^1_f = 0$. We see that from $(\ddagger)$, we must have $\big(|\Sigma f|\cdot V(z_0)\big)_\bold 0 =1$ in this case.

\vskip .2in

\noindent$\bullet$\hskip .2in Consider now the next easiest case where $\mu_0-\lambda^1_f = 1$. Then, $(\ddagger)$ implies that $\big(|\Sigma f|\cdot V(z_0)\big)_\bold 0$ is equal to $1$ or $2$. However, $(\dagger)$ tells us that we cannot have $\big(|\Sigma f|\cdot V(z_0)\big)_\bold 0=1$. 

Therefore, we must have $\big(|\Sigma f|\cdot V(z_0)\big)_\bold 0=2$. This implies that $\Sigma f$ has two smooth components which are each transversely intersected by $V(z_0)$.

\vskip .2in

By taking the contrapositive of the above work, we recover a slightly-improved version of the result of [{\bf Ma1}]: if  $\mu_0-\lambda^1_f = 1$, then either $\Sigma f$ consists of two smooth components which are transversely intersected by $V(z_0)$ at the origin, or the rank of $\widetilde H^{n-1}(F; \ \Bbb Z)$ is strictly less than $\lambda^1_f$.

\vskip .3in

\noindent{\bf Application 3}. Suppose that $f_0$ is homogeneous of degree $d_0$. By using Theorem 9.6 of [{\bf Mi}], or directly from [{\bf M-O}], the characteristic polynomial of the $f_0$-monodromy is given by
$$
{\operatorname{char}}_{h_0}(t) = (t-1)^{a_0}\left(\frac{t^{d_0}-1}{t-1}\right)^{b_0},\tag{$\dagger$}
$$
where $\dsize b_0:= \frac{(d_0-1)^n-(-1)^n}{d_0}$ and $\dsize a_0:= b_0+(-1)^n = (d_0-1)\left(\frac{(d_0-1)^{n-1}-(-1)^{n-1}}{d_0}\right)$.

\vskip .1in

Let $\Phi_k$ denote the $k$-th cyclotomic polynomial, i.e., $\Phi_k=\prod_{\xi}(t-\xi)$ where the $\xi$ vary over the primitive $k$-th roots of unity. Then, the main theorem immediately tells us that 
$$
{\operatorname{char}}_{m_{n-1}}(t) = \prod_{k|d_0} \Phi_k^{c_k},
$$
where $c_1\leqslant a_0$ and, for $k>1$, $\dsize c_k\leqslant b_0$; in addition, for all $k|d_0$, we must have that $\Phi_k^{c_k}$ divides $\prod_\nu{\operatorname{char}}_{h_\nu}(t)$ in $\Bbb Z[t]$.

\vskip .1in

This is particularly useful in the special case where, at each point $\bold p_i\in H_t\cap\Sigma f$, $f_t$ is ``homogeneous at $\bold p_i$'', i.e., $g_i(\bold z):=f_t(\bold z+\bold p_i)$ is homogeneous. Let $d_i$ denote the degree of $g_i$. For all $\bold p_i$ in a given component $\nu$ of $\Sigma f$, the $d_i$ must be the same; denote this common value by $d_\nu$.
For each $\nu$, there is a formula analogous to $(\dagger)$ for the characteristic polynomial of $h_\nu$, and so the irreducible factors of $\prod_\nu{\operatorname{char}}_{h_\nu}(t)$ are cyclotomic polynomials $\Phi_k$ for which $k$ must divide one of the $d_\nu$.

An example of such an $f$ would be one which defines a central arrangement of $d_0$ hyperplanes in $\Bbb C^3$. 
The above paragraphs generalize most of the main theorem from [{\bf Ma3}], where we considered only hyperplane arrangements in $\Bbb C^3$. This is related to the work of [{\bf E}].

\vskip .3in

\noindent\S4. {\bf Further Remarks}.

\vskip .1in
Regardless of the dimension of the critical locus of $f$, a nexus diagram exists in the category of perverse sheaves. 

In the proof of the Main Theorem, even before we took the stalk cohomology at the origin, we had two short exact sequences in the category of perverse sheaves; that is, we had a perverse nexus diagram. However, since each of these perverse sheaves had the origin as an isolated point in their support, taking stalk cohomology essentially did nothing.

If the dimension of $\Sigma f$ is arbitrary and $z_0$ is generic, we still obtain a perverse nexus diagram 
$$
\CD
{}@. \psi_{\hat z_0}[-1]\phi_f[-1]\Bbb Z^\bullet_{\Cal U}[n+1] @.{} \\
@. @VV\alpha V @. \\
\phi_{f_0}[-1]\Bbb Z_{\Cal U\cap V(z_0)}^\bullet[n]@>\gamma>>\bold N^\bullet@>\delta>>j_*\Bbb Z_{\bold 0}^{\lambda^0_{f}}\\
@. @VV\beta V @. \\
{}@. j_*\Bbb Z_{\bold 0}^{\omega}, @.{}
\endCD
$$
where $j$ denotes the inclusion of the origin into $V(f, z_0)$. The $f$ and $z_0$ monodromies act commutatively on this diagram.

We write ${\operatorname{ker}}^p$ and ${\operatorname{coker}}^p$ for the kernel and cokernel in the category of perverse sheaves, and we write ${}^p\hskip -.02in H^*$ for the perverse cohomology. Then, in general, in the perverse nexus diagram above,
$$
{\operatorname{ker}}^p(\beta\circ\gamma)\ \cong\ {\operatorname{ker}}^p(\delta\circ\alpha)\ \cong\ {}^p\hskip -.02in H^{-1}\big(\hat k^*\phi_f[-1]\Bbb Z^\bullet_{\Cal U}[n+1]\big)
$$
and 
$$
{\operatorname{coker}}^p(\beta\circ\gamma)\ \cong\ {\operatorname{coker}}^p(\delta\circ\alpha)\ \cong\ {}^p\hskip -.02in H^{0}\big(\hat k^*\phi_f[-1]\Bbb Z^\bullet_{\Cal U}[n+1]\big).
$$

\vskip .2in

Furthermore, there is an analogous perverse nexus diagram involving $\psi_f[-1]$, in place of $\phi_f[-1]$. That is, there is a nexus diagram in the category of perverse sheaves obtained from the above diagram by replacing each occurrence of  $\phi_f[-1]$ by $\psi_f[-1]$, and by replacing $\lambda^0_f$ by $\omega$. The perverse kernel and cokernel statements also hold with $\phi_f[-1]$ replaced by $\psi_f[-1]$.

\enddocument